\documentclass[12pt]{amsart}
\usepackage{amscd}      
\usepackage{amssymb}
\usepackage{amsmath}
\usepackage{xypic}      
\LaTeXdiagrams          
\usepackage[all,v2]{xy}
\xyoption{2cell} \UseAllTwocells \xyoption{frame} \CompileMatrices
\allowdisplaybreaks[3]

\usepackage{latexsym}
\usepackage{epsfig}
\usepackage{amsfonts}
\usepackage{enumerate}

\newtheorem{prop}{Proposition}[section]

\newtheorem{thm}[prop]{Theorem}

\newtheorem{conjecture}[prop]{Conjecture}

\newtheorem{rmk}{Remark}

            {\nolinebreak $\Box$ \end{trivlist}}

\newcommand{\noprint}[1]{}

\newcommand{\zz}{{\mathbb Z}}

\newcommand{\cc}{{\mathbb C}}

\newcommand{\Aut}{\mathop{\rm Aut}\nolimits}

\newcommand{\ldiag}[1]%
       {\makebox[0cm]{${\scriptstyle#1}\downarrow\phantom{\scriptstyle#1}$}}
\newcommand{\ldiagup}[1]%
       {\makebox[0cm]{${\scriptstyle#1}\uparrow\phantom{\scriptstyle#1}$}}
\newcommand{\rdiag}[1]%
       {\makebox[0cm]{$\phantom{\scriptstyle#1}\downarrow{\scriptstyle#1}$}}
\newcommand{\sediagr}[1]%
       {\makebox[0cm]{$\phantom{\scriptstyle#1}\searrow{\scriptstyle#1}$}}
\newcommand{\nediagr}[1]%
       {\makebox[0cm]{$\phantom{\scriptstyle#1}\nearrow{\scriptstyle#1}$}}
\newcommand{\rdiagup}[1]%
       {\makebox[0cm]{$\phantom{\scriptstyle#1}\uparrow{\scriptstyle#1}$}}
\newcommand{\swdiag}[1]%
       {\makebox[0cm]{$\phantom{\scriptstyle#1}\swarrow{\scriptstyle#1}$}}
\newcommand{\sediag}[1]%
       {\makebox[0cm]{${\scriptstyle#1}\searrow\phantom{\scriptstyle#1}$}}
\newcommand{\nediag}[1]%
       {\makebox[0cm]{${\scriptstyle#1}\nearrow\phantom{\scriptstyle#1}$}}

\newcommand{\doublearrowstack}[2]%
                      {{{{\scriptstyle#1}\atop{\textstyle\longrightarrow}}\atop{{\textstyle\longrightarrow}\atop{\scriptstyle#2}}}}
\newcommand{\rightleftarrowstack}[2]%
                      {{{{\scriptstyle#1}\atop{\textstyle\longrightarrow}}\atop{{\textstyle\longleftarrow}\atop{\scriptstyle#2}}}}
\newcommand{\leftrightarrowstack}[2]%
                      {{{{\scriptstyle#1}\atop{\textstyle\longleftarrow}}\atop{{\textstyle\longrightarrow}\atop{\scriptstyle#2}}}}

\newcommand{\overtoparrow}%
{\makebox[0cm]{\beginpicture \setcoordinatesystem units
<.8cm,.4cm> point at 0 0 \setplotarea x from -3 to 3, y from 0 to
1 \setquadratic \plot -3 0 0 1 3 0 / \put{\vector(3,-1){0}}[Bl] at
3 0
\endpicture}}

\newcommand{\underbottomarrow}%
{\makebox[0cm]{\beginpicture \setcoordinatesystem units
<.8cm,.4cm> point at 0 0 \setplotarea x from -3 to 3, y from 0 to
1 \setquadratic \plot -3 1 0 0 3 1 / \put{\vector(3,1){0}}[Bl] at
3 1
\endpicture}}

\newcommand{\ses}[5]%
{0\longrightarrow#1\stackrel{#2}{ \longrightarrow}#3\stackrel{#4}{
\longrightarrow}#5\longrightarrow0}

\newcommand{\dt}[6]%
{#1\stackrel{#2}{longrightarrow}#3
\stackrel{#4}{\longrightarrow}#5 \stackrel{#6}{\longrightarrow}
#1[1]}

\newcommand{\cat}[1]%
{(\mbox{\rm #1})}

\def\Label#1{\label{#1}{\tt [#1]}\phantom{h}}
\def\Label{\label}

\newtheorem{lemma}[subsubsection]{Lemma}

\theoremstyle{definition}

\newtheorem{definition}[subsubsection]{Definition}

\theoremstyle{remark}

\theoremstyle{remark}

\numberwithin{equation}{section}

\newcommand{\Mbar}{\overline{\M}}

\newcommand{\com}{\mathbb{C}}

\newcommand{\X}{\mathcal{X}}

\newcommand{\M}{\mathcal{M}}
\newcommand{\C}{\mathcal{C}}

\newcommand{\B}{\mathcal{B}}

\newcommand{\D}{\mathcal{D}}
\newcommand{\F}{\mathcal{F}}

\newcommand{\sL}{\mathcal{L}}
\newcommand{\sI}{\mathcal{I}}

\newcommand{\bt}{\mathbf{t}}
\newcommand{\bq}{\mathbf{q}}
\newcommand{\sH}{\mathcal{H}}

\newcommand{\bp}{\mathbf{p}}

\def\<{\left\langle}
\def\>{\right\rangle}

\title{On Virasoro Constraints for Orbifold Gromov-Witten Theory}
\author{Yunfeng Jiang}
\address{Department of Mathematics\\ University of British Columbia\\ 1984 Mathematics Road\\Vancouver\\ BC V6T 1Z2\\ Canada}
\email{jiangyf@math.ubc.ca}
\author{Hsian-Hua Tseng}
\address{Department of Mathematics\\ University of British Columbia\\ 1984 Mathematics Road\\Vancouver\\ BC V6T 1Z2\\ Canada}
\email{hhtseng@math.ubc.ca}
\date{\today}

\begin{document}
\begin{abstract}
Virasoro constraints for  orbifold Gromov-Witten
theory are described. These constraints are applied to the degree zreo, genus zero
orbifold Gromov-Witten potentials of the weighted projective stacks $\mathbb{P}(1,N)$, $\mathbb{P}(1,1,N)$ and $\mathbb{P}(1,1,1,N)$ to obtain formulas of descendant cyclic Hurwitz-Hodge integrals.
\end{abstract}

\maketitle

\section{Introduction}
One of the interesting conjectures in Gromov-Witten theory is the so-called Virasoro conjecture, proposed by T. Eguchi, K. Hori, M. Jinzenji, C.-S. Xiong, and S. Katz in \cite{ehx, ejx}. There has been significant developments in understanding this conjecture. To the best of our knowledge, Virasoro conjecture has been proven for toric manifolds (\cite{gi_gwi}, \cite{Ir}), flag manifolds \cite{KJ}, Grassmanians \cite{BCFK}, and nonsingular curves \cite{OP_vir}. Recent breakthrough by C. Teleman should soon lead to a proof of Virasoro conjecture for compact K\"ahler manifolds with semi-simple quantum cohomologies.

B. Dubrovin and Y. Zhang \cite{DZ} have defined Virasoro operators associated to a general Frobenius manifold. This can be viewed as promoting Virasoro constraints to an axiom of the abstract topological field theory.   From this point of view it is not surprising that Virasoro constraints can be formulated for Gromov-Witten theory of K\"ahler {\em orbifolds}. Our purpose of this paper is to explicitly describe these conjectural constraints, and apply them to study degree zero invariants, along the line of \cite{GP}.

The rest of this paper is organized as follows. Section \ref{review} contains a very brief recollection of ingredients from orbifold Gromov-Witten theory needed in the paper. In Section \ref{virasoro} we describe the conjectural Virasoro constraints for Gromov-Witten theory of K\"ahler orbifolds, and discuss some related issues. Sections \ref{degree_0_maps} through \ref{threefolds} are devoted to the applications of Virasoro constraints to genus zero degree zero invariants. In Appendix \ref{conj_orbLW} we discuss a conjecture concerning certain genus one invariants and characteristic numbers.

\subsection*{Acknowledgments}
 H.-H. T. is grateful to A. Givental for his generous help on his approach to Virasoro constraints, and to E. Getzler for suggesting the problem of applying Virasoro constraints to degree zero invariants in May, 2004. Y. J. was supported by the University Graduate Fellowship from the University of British Columbia. H.- H. T. is supported in part by a postdoctoral fellowship from the Pacific Institute of Mathematical Sciences (Vancouver, Canada) and a visiting research fellowship from Institut Mittag-Leffler (Djursholm, Sweden). 

\section{Preliminaries}\label{review}
In this section we review some materials essential to the construction of orbifold Gromov-Witten theory. Gromov-Witten theory for symplectic orbifolds were constructed in \cite{CR1, CR2} and for Deligne-Mumford stacks in \cite{AGV, agv2}. By now there exist several articles \cite{Tseng}, \cite{cclt}, \cite{ccit} containing summaries of basic materials. So we should refrain ourselves from repetition and only give a very brief summary of the ingredients we need. We refer to \cite{Tseng}, \cite{cclt}, \cite{ccit} , \cite{CR1, CR2}, \cite{AGV, agv2} for detail discussions.

Throughout this paper let $\X$ denote a compact K\"ahler orbifold. The main ingredients of orbifold theory are summarized as follows:

\begin{center}
  \begin{tabular}{lp{0.75\textwidth}}
    $I\X$ & the inertia orbifold of $\X$.  A point of $I\X$ is a pair $(x,g)$ with $x$ a point of $\X$ and
    $g\in \Aut_{\X}(x)$.\\
    $I\X = \coprod_{i \in \sI} \X_i$ & the decomposition of $I\X$
    into components; here $\sI$ is an index set. \\
    $q$ & the natural projection $I\X\to \X$ defined by forgetting $g$.\\
    $I$ & the involution of $I\X$ which sends $(x,g)$ to $(x,g^{-1})$. \\
    $i^I$ & the index associated to $i\in \sI$ so that $\X_i$ and $\X_{i^I}$ are isomorphic under the involution $I$.\\
    $H_{CR}^*(\X,\com)$ & the Chen-Ruan orbifold cohomology groups of $\X$.  These are the
    cohomology groups $H^*(I\X;\com)$ of the inertia orbifold.  \\
    $age$ & a rational number associated to each component $\X_i$ of
    the inertia stack. This is called the degree-shifting number in \cite{CR1, CR2}. \\
    $orbdeg$ & the orbifold degree: For a class $a\in H^p(\X_i)$,  define $orbdeg (a):=p+2age(\X_i)$. \\
    $(\alpha,\beta)_{orb}$ & the
    orbifold Poincar\'e pairing $\int_{I\X} \alpha \cup I^* \beta$. \\
    $\{\phi_\alpha\}$ & an additive homogeneous basis of $H_{CR}^*(\X,\com)$.\\
    $\phi_0=1$ & the class in $H^0(\X,\com)$ Poincar\'e dual to the fundamental class.\\
    $\{\phi^\alpha\}$ & the basis of $H_{CR}^*(\X,\com)$ dual under orbifold Poincar\'e pairing.\\

  \end{tabular}
\end{center}

Here is a summary of the main ingredients of orbifold Gromov-Witten theory:

\begin{center}
 \begin{tabular}{lp{0.75\textwidth}}
$\Mbar_{g,n}(\X,\beta)$ & the moduli stack of genus $g$, $n$-pointed orbifold stable maps of degree $\beta$, with sections to all gerbes.  \\
$ev_j$ & the evaluation map $\Mbar_{g,n}(\X,\beta)\to I\X$ associated to the $j$-th marked point.\\
$\Mbar_{g,n}(\X,\beta,(i_1,..,i_n))$ & the open-and-closed substack $ev_1^{-1}(\X_{i_1})\cap...\cap ev_n^{-1}(\X_{i_n})$.\\
$[\Mbar_{g,n}(\X,\beta, (i_1,...,i_n))]^{w}$ & the (weighted) virtual fundamental class.\\
$\psi_j$ & the descendent class.\\
$\<a_1\psi_1^{k_1}...a_n\psi_n^{k_n}\>^\X_{g,n,\beta}$  & the descendent orbifold Gromov-Witten invariant, which is defined to be $\int_{[\Mbar_{g,n}(\X,\beta,(i_1,...,i_n))]^{vir}}ev_1^*a_1\psi_1^{k_1}\wedge...\wedge ev_n^*a_n\psi_n^{k_n}$.\\
$*$ & the orbifold cup product of Chen-Ruan.\\
$\Lambda$ & the Novikov ring associated to $H_2(X,\mathbb{Z})$.\\
$\star_t$ & the orbifold quantum product, depending on $t\in H_{CR}^*(\X,\com)$.\\
$t_k^\alpha$ & coordinates of cohomology under the basis $\{\phi_\alpha\}$: $t_k=\sum_\alpha t_k^\alpha\phi_\alpha$.\\
$\bt$ & $\sum_{k\geq 0} t_k z^k$.\\

 \end{tabular}
\end{center}
The {\em total descendent potential} is defined to be $$\D_\X(\bt):=\exp\left(\sum_{g\geq 0} \hbar^{g-1}\F_\X^g(\bt)\right),$$ where
$$\F_\X^g(\bt):=\sum_{n,d}\frac{Q^{d}}{n!}\<\mathbf{t},...,\mathbf{t}\>_{g,n,d}=\sum_{n,d}\frac{Q^d}{n!}\int_{[\Mbar_{g,n}(\X,d)]^{vir}}\bigwedge_{i=1}^n \sum_{k=0}^\infty ev_i^*t_k\psi_i^k.$$
$\F_\X^g(\bt)$ is called the genus-$g$ descendent potential, it is well-defined as a $\Lambda$-valued formal power series in coordinates $t_k^\alpha$. The total descendent potential $\D_\X(\bt)$ is well-defined as a formal power series in $t_k^\alpha$ taking values in $\Lambda[[\hbar,\hbar^{-1}]]$.

\section{Virasoro Constraints for Orbifolds}\label{virasoro}
In this Section we formulate the conjectural Virasoro constraints for orbifold Gromov-Witten theory. We follow the approach of \cite{DZ} and \cite{gi_gwi}.

\subsection{Quantization Formalism}\label{giv}
The material of this section is due to Givental, and adapted to orbifold Gromov-Witten theory in \cite{Tseng}.

Consider the space $$\sH:=H^*(I\X, \com)\otimes\Lambda\{z,z^{-1}\}$$ of orbifold-cohomology-valued convergent Laurent series (see \cite{ccit2}, Section 3). There is a symplectic form on $\sH$ given by
$$\Omega(f,g)=\text{Res}_{z=0} (f(-z),g(z))_{orb}dz, \quad \text{for } f,g\in\sH.$$
There is a polarization $$\sH_+:= H^*(I\X,\com)\otimes \Lambda\{z\},\quad \sH_-:=z^{-1}H^*(I\X,\com)\otimes \Lambda\{z^{-1}\}.$$

This identifies $\sH$ with the cotangent bundle $T^*\sH_+$. Both $\sH_+$ and $\sH_-$ are Lagrangian subspaces with respect to $\Omega$.

Introduce the Darboux coordinate system on $(\sH,\Omega)$ with respect to the polarization above, $$\{p_a^\mu,q_b^\nu\}.$$
In these coordinates, a general point in $\sH$ takes the form
$$\sum_{a\geq 0}\sum_\mu p_a^\mu\phi^\mu(-z)^{-a-1}+\sum_{b\geq 0}\sum_\nu q_b^\nu\phi_\nu z^b.$$
Put $p_a=\sum_\mu p_a^\mu\phi^\mu$ and $q_b=\sum_\nu q_b^\nu \phi_\nu$. Denote $$\bp=\bp(z):=\sum_{k\geq 0} p_kz^{-k-1}=p_0z^{-1}+p_1z^{-2}+..., \quad \text{and } \bq=\bq(z):=\sum_{k\geq 0}q_kz^k=q_0+q_1z+q_2z^2+....$$

For $\bt(z)\in\sH_+$ introduce a shift $\bq(z)=\bt(z)-1z$ called the dilaton shift. Define the Fock space $\mbox{\it Fock}$ to be the space of formal functions in $\bq(z)=\bt(z)-1z\in \sH_+$. In other words, $\mbox{\it Fock}$ is the space of formal functions on $\sH_+$ near $\bq=-1z$. The descendent potential $\D_\X(\bt)$ is regarded as an {\em asymptotical element} (see \cite{gi_An}, Section 8 for the definition) in {\it Fock} via the dilaton shift.

Let $A:\sH\to\sH$ be a linear infinitesimally symplectic transformation, i.e. $\Omega (Af, g)+\Omega (f, Ag)=0$ for all $f, g\in \sH$. When $f\in \sH$ is written in Darboux coordinates, the quadratic Hamiltonian $$f\mapsto\frac{1}{2}\Omega(Af,f),$$ is a series of homogeneous degree two monomials in Darboux coordinates $p_a^\alpha,q_b^\alpha$ in which each variable occurs only finitely many times. Define the quantization of quadratic monomials as
$$\widehat{q_a^\mu q_b^\nu}=\frac{q_a^\mu q_b^\nu}{\hbar},\,\, \widehat{q_a^\mu p_b^\nu}=q_a^\mu\frac{\partial}{\partial q_b^\nu},\,\,\widehat{p_a^\mu p_b^\nu}=\hbar\frac{\partial}{\partial q_a^\mu}\frac{\partial}{\partial q_b^\nu}.$$
Extending linearly, this defines a quadratic differential operator $\widehat{A}$ on $\mbox{\it Fock}$, called the quantization of $A$. The differential operators $\widehat{q_aq_b},\widehat{q_ap_b}, \widehat{p_ap_b} $ act on $\mbox{\it Fock}$. Since the quadratic Hamiltonian of $A$ may contain infinitely many monomials, the quantization $\widehat{A}$ does not act on ${\it Fock}$ in general. 

For infinitesimal symplectomorphisms $A$ and $B$, there is the following relation
$$[\widehat{A},\widehat{B}]=\{ A,B \}^\wedge+ \C (h_A,h_B),$$ where $\{\cdot,\cdot\}$ is the Lie bracket, $[\cdot,\cdot]$ is the supercommutator, and $h_A$ (respectively $h_B$) is the quadratic Hamiltonian of $A$ (respectively $B$). A direct calculation shows that the cocycle $\C$ is given by
\begin{equation*}
\begin{split}
&\C(p_a^\mu p_b^\nu,q_a^\mu q_b^\nu)=-\C(q_a^\mu q_b^\nu,p_a^\mu p_b^\nu)=1+\delta^{\mu\nu}\delta_{ab},\\
&\C=0\mbox{ on any other pair of quadratic Darboux monomials.}
\end{split}
\end{equation*}
For simplicity, we write $\C(A,B)$ for $\C (h_A,h_B)$.

\subsection{Virasoro Constraints}
As explained in e.g. Section 3.2 of \cite{CR1}, cohomology groups $H^*(\X_i)$ of components of $I\X$ admits Hodge decompositions $$H^k(\X_i, \com)=\oplus_{p+q=k}H^{p,q}(\X_i,\com).$$

From now on, the additive basis $\{\phi_\alpha\}$ will be assumed homogeneous with respective to the Hodge decomposition.

Define two operators $\rho$ and $\mu$ as follows:
\begin{enumerate}
\item
$\rho:=c_1(T_\X)* : H^*(I\X,\Lambda)\to H^*(I\X,\Lambda)$ is defined to be the orbifold multiplication by the first Chern class $c_1(T_\X)$.
\item
$\mu: H^*(I\X,\Lambda)\to H^*(I\X,\Lambda)$ is defined as follows: for a class $\alpha\in H^{p,q}(\X_i,\com)$, define $$\mu(\alpha):=\left(p+age(\X_i)-\frac{\text{dim}\X}{2}\right)\alpha.$$
\end{enumerate}
We write $\rho_{\alpha}^\beta$ for the matrix of $\rho$ under the basis $\{\phi_\alpha\}$:
$$\sum_\beta\rho_\alpha^\beta\phi_\beta= \rho(\phi_\alpha).$$
In the basis $\{\phi_\alpha\}$, the operator $\mu$ is given by a diagonal matrix with entries $\mu_\alpha$: $$\mu(\phi_\alpha)=\mu_\alpha\phi_\alpha.$$

The following is easy to check:
\begin{lemma}
\hfill
\begin{enumerate}
\item
$(\mu(a), b)_{orb}+(a,\mu(b))_{orb}=0$;
\item
$(\rho(a),b)_{orb}=(a,\rho(b))_{orb}$.
\end{enumerate}
\end{lemma}
\begin{proof}
(2) is obvious. To check (1), we may assume that $a\in H^{p,q}(\X_i,\com)$ and $b\in H^{\text{dim}\X_i-p, \text{dim}\X_i-q}(\X_{i^I},\com)$. Then we have
\begin{equation*}
\begin{split}
&(\mu(a),b)_{orb}+(a,\mu(b))_{orb}\\
&=\left(p+age(\X_i)-\frac{\text{dim}\X}{2}+\text{dim}\X_i-p+age(\X_{i^I})-\frac{\text{dim}\X}{2}\right)(a,b)_{orb}\\
&=\left(age(\X_i)+age(\X_{i^I})+\text{dim}\X_i-\text{dim}\X\right)(a,b)_{orb}=0,
\end{split}
\end{equation*}
where the last equality follows from \cite{CR1}, Lemma 3.2.1.
\end{proof}

In Givental's approach to Virasoro constraints \cite{gi_gwi}, the Virasoro operators are written using the quantization formalism reviewed in Section \ref{giv}. This is what we do next.

\begin{definition}
For $m\geq -1$, put $$L_m^{\mu,\rho}:=z^{-1/2}\left(z\frac{d}{dz}z-\mu z+\rho\right)^{m+1}z^{-1/2}, $$
and define
$$\sL_m:=(L_m^{\mu,\rho})^\wedge+\frac{\delta_{m,0}}{4}str\left(\frac{1}{4}-\mu\mu^*\right).$$
\end{definition}

\begin{lemma}
The operators $\sL_m$ satisfy the commutation relations $$[\sL_m,\sL_n]=(m-n)\sL_{m+n}.$$
\end{lemma}

\begin{proof}
It is clear that the operators $L_m^{\mu,\rho}$ satisfy $$[L_m^{\mu,\rho},L_n^{\mu,\rho}]=(m-n)L_{m+n}^{\mu, \rho}.$$ Therefore the quantized operators $(L_m^{\mu,\rho})^\wedge$, corrected by some central constants $c_m$, satisfy the same commutation relations. A direct computation shows that the constants are $c_m=0$ for $m\neq 0$ and $c_0=\frac{\delta_{m,0}}{4}str(\frac{1}{4}-\mu\mu^*)$.
\end{proof}

It is clear that when $\X$ is a K\"ahler manifold, $\sL_m$ coincide with the Virasoro operators written in \cite{gi_gwi}. Indeed, our construction is a straightforward adaptation of that in \cite{gi_gwi} and Dubrovin-Zhang's construction \cite{DZ} for Frobenius manifolds.

One can write down the operators $\sL_m$ explicitly. Let $r:=\text{dim}\X$. Define the symbol $[x]_i^k$ by $$\sum_{i=0}^{k+1} s^i[x]_i^k=\prod_{i=0}^k(s+x+i).$$
Then we have:
\begin{multline*}
\sL_m = \sum_{i=0}^{m+1}\left(\frac{\hbar}{2}\sum_{k=i-m}^{-1}(-1)^k[\mu_\alpha+k+\frac{1}{2}]_i^m(\rho^i)^{\alpha\beta}\frac{\partial}{\partial t_{-k-1}^\alpha}\frac{\partial}{\partial t_{k+m-i}^\beta}
- [\frac{3-r}{2}]_i^m(R^i)_0^b\frac{\partial}{\partial t_{m-i+1}^\beta}\right. \\
+ \left.\sum_{k=0}^\infty [\mu_\alpha+k+\frac{1}{2}]_i^m(R^i)_\alpha^\beta t_k^\alpha\frac{\partial}{\partial t_{k+m-i}^\beta}\right)+\frac{1}{2\hbar}(R^{m+1})_{\alpha\beta}t_0^\alpha t_0^\beta+\frac{\delta_{m,0}}{4}str\left(\frac{1}{4}-\mu\mu^*\right).
\end{multline*}
This can be seen as follows. Note that the right side is exactly what comes out of Dubrovin-Zhang's construction. Therefore both sides satisfy the Virasoro commutation relations. For this reason, it suffices to check the equation above for $\sL_{-1}$ and $\sL_2$, which can be done by a direct computation.

\begin{conjecture}[Virasoro constraints]
\begin{equation}\label{vir_constraint}
\sL_m\D_\X=0, \text{ for } m\geq -1.
\end{equation}
\end{conjecture}

\subsection{Relation to Geometric Theory}
The Virasoro operators $\sL_m$ are constructed without reference to orbifold Gromov-Witten theory. In this section we discuss the compatibility of (\ref{vir_constraint}) with what's known from geometric theory.

Observe that the operator $\sL_{-1}\D_\X=\widehat{(1/z)}\D_\X=0$ is the {\em string equation}, which always holds. See \cite{Tseng} for more discussion on the string equation in orbifold Gromov-Witten theory.

In Gromov-Witten theory of K\"ahler manifolds, the constraint $\sL_0\D_\X=0$, known as Hori's equation, is derived by combining the virtual dimension formula, the divisor equation, and the dilaton equation. The same argument implies

\begin{lemma}[c.f. \cite{G}, Theorem 2.1]
We have $\sL_0'\D_\X=0$, where
\begin{multline*}
\sL_0'=-\frac{1}{2}(3-r)\frac{\partial}{\partial t_1^{0}}+\sum_{m=0}^\infty(\mu_\alpha+m+\frac{1}{2})t_m^\alpha \frac{\partial}{\partial t_m^\alpha}-R_0^\beta\frac{\partial}{\partial t_0^\beta} \\
+\sum_{m=1}^\infty R_\alpha^\beta t_m^\alpha\frac{\partial}{\partial t_{m-1}^\beta}+\frac{1}{2\hbar}R_{\alpha\beta}t_0^\alpha t_0^\beta+ \frac{1}{2}(3-r)\<\psi\>_{1,1,0}^\X-\<c_1(T_\X)\>_{1,1,0}^\X.
\end{multline*}
\end{lemma}
\begin{proof}
This is done exactly as the proof of \cite{G}, Theorem 2.1. We simply remark that the terms $\<\psi\>_{1,1,0}^\X$ and $\<c_1(T_\X)\>_{1,1,0}^\X$ come from the exceptional cases of the dilaton and divisor equations respectively.
\end{proof}

Observe that the non-constant parts of the operators $\sL_0$ and $\sL_0'$ are the same. Thus the consistency with geometric theory forces the constant terms to be the same. This imposes the following

\begin{conjecture}\label{orbLW}
$$\frac{1}{4}str\left(\frac{1}{4}-\mu^2\right)=\frac{1}{2}(3-r)\<\psi\>_{1,1,0}^\X-\<c_1(T_\X)\>_{1,1,0}^\X.
$$
\end{conjecture}
For evidences and discussions of Conjecture \ref{orbLW}, see Appendix \ref{conj_orbLW}.

\subsection{Evidence}
In this section we discuss some evidence of the Virasoro constraints (\ref{vir_constraint}).

The first evidence is the case $\X=\mathcal{B}G$ for a finite group $G$. The orbifold Gromov-Witten theory of $\mathcal{B}G$ has been completely solved by Jarvis-Kimura \cite{jk}. Among other things, the Virasoro constraints for $\mathcal{B}G$ is explicitly written down and proven. It is clear that the constraints in \cite{jk} coincide with (\ref{vir_constraint}). In particular, Conjecture \ref{orbLW} is verified for $\mathcal{B}G$ by evaluating the right side.

The second evidence is that Virasoro constraints hold in genus zero. It is known (e.g. \cite{ehx}, \cite{G}) that Virasoro constraints in genus zero Gromov-Witten theory is a formal consequence of the string equation, dilaton equation, and topological recursion equations. A geometric proof based on Givental's symplectic space approach to Gromov-Witten theory is given in \cite{gi_frob}, Theorem 6 for K\"ahler manifolds. The same argument proves the case of K\"ahler orbifolds as well. This is sufficient for the applications of Virasoro constraints discussed in later sections.

Finally, the constraints (\ref{vir_constraint}) will be proven in \cite{mt} for weighted projective lines $\mathbb{P}(k,m)$ for $k,m$ co-prime.

\section{Degree zero twisted stable maps to orbifolds.}\label{degree_0_maps}

In this section we discuss the degree zero twisted stable maps to
orbifolds. We describe the degree zero orbifold Gromov-Witten
invariants as Hurwitz-Hodge integrals. Let $\Sigma_{i}n_{i}:=n_{1}+\cdots+n_{N-1}$, where
$n_{1},\cdots,n_{N-1}$ are positive integers. Throughout the paper we use
the notation $k_{1}\cdots k_{\Sigma_{i}n_{i}}$ to represent $k_{1}\cdots k_{n_{1}}k_{n_{1}+1}\cdots k_{n_{2}}k_{n_{2}+1}\cdots k_{\Sigma_{i}n_{i}}$.
For example, $\widetilde{\tau}_{k_{1}}\cdots\widetilde{\tau}_{k_{\Sigma_{i}
n_{i}}}$
represents
$\widetilde{\tau}_{k_{1}}\cdots\widetilde{\tau}_{k_{n_{1}}}\widetilde{\tau}_{k_{n_{1}+1}}\cdots\widetilde{\tau}_{k_{n_{2}}}
\widetilde{\tau}_{k_{n_{2}+1}}\cdots\widetilde{\tau}_{k_{\Sigma_{i}
n_{i}}}$.

Let $\mathcal{X}=\mathbb{P}(1,\cdots,1,N)$ be the weighted projective stack of dimension $d$ with weights $(1,1,\cdots,1,N)$. It is an orbifold which has only one orbifold point $p=[0,\cdots,0,1]$ with local group the cyclic group $\mathbb{Z}_{N}$. We consider the degree zero twisted stable maps to $\mathcal{X}$. From now on, put $\omega=\exp(\frac{2\pi\sqrt{-1}}{N})$.

Let $f: \mathcal{C}\rightarrow \mathcal{X}$ be a degree zero twisted
stable map. If there are stacky points on $\C$, then $f$ factors through $p\simeq \mathcal{B}\zz_N$. Suppose that there are $n_{1}+n_{2}+\cdots+n_{N-1}$ stacky points on the curve $\mathcal{C}$.
The orbifold fundamental group of $\C$ can be presented as
$$\pi^{orb}_{1}(\mathcal{C})=\{\xi_{1},\cdots,\xi_{\Sigma_{i}n_{i}}: \xi_{j}^{a_{j}}=1, \xi_{1}\cdots\xi_{\Sigma_{i}n_{i}}=1\},$$
where $\xi_{j}$ is a generator represented by a loop around the $j$-th stacky point and $a_j$
is the local index of the $j$-th stacky point. The map $f$ induces a homomorphism
$$\varphi: \pi^{orb}_{1}(\mathcal{C})\to \mathbb{Z}_{N}$$ 
which is given by $\xi_{j}\mapsto \omega^{b_{j}}$  for $1\leq j\leq\Sigma_{i}n_{i}$ satisfying the condition that $0<b_{j}\leq N-1$ and $a_{j}b_{j}$ is divisible by $N$. 
The morphism $f$ is equivalent to having an admissible $\mathbb{Z}_{N}$-cover $\widetilde{\mathcal{C}}\to C$ over the {\em coarse curve} with ramification condition specified by $\varphi$. Suppose the 
stacky type on the curve $\mathcal{C}$ is given  as follows,
\begin{equation}\label{ramification}
\mathbf{x}=\left(\underbrace{\omega,\cdots,\omega}_{n_{1}},\underbrace{\omega^{2},\cdots,\omega^{2}}_{n_{2}},
\cdots,\underbrace{\omega^{i},\cdots,\omega^{i}}_{n_{i}},\cdots,\underbrace{\omega^{N-1},\cdots,\omega^{N-1}}_{n_{N-1}}\right).
\end{equation}
This means that $\varphi(\xi_j)=\omega^i$ for $\sum_{a=0}^{i-1}n_a+1\leq j\leq \sum_{a=0}^i n_a$ (we put $n_0=0$).

The equation $\prod_j\xi_j=1$ in $\pi_1^{orb}(\C)$ translates into the following condition on $n_{1},\cdots,n_{N-1}$:
$$n_{1}+2n_{2}+\cdots+(N-1)n_{N-1}\equiv 0\,(\text{mod }N).$$

Let $\overline{\mathcal{M}}_{n+\Sigma_{i}n_{i}}(\mathcal{X},0,\mathbf{x})$ be
the moduli stack of degree zero, genus zero twisted stable maps with
$n$ non-stacky marked points and $\Sigma_{i}n_{i}$ stacky marked points
to $\mathcal{X}$ such that the stacky structure on the orbicurve $\mathcal{C}$ is given by $\mathbf{x}$ in (\ref{ramification}). Let 
$$\pi: \overline{\mathcal{M}}_{n+\Sigma_{i}n_{i}}(\mathcal{X},0,\mathbf{x})\to \overline{M}_{0,n+\Sigma_{i}n_{i}}$$
be the forgetful map, where $\overline{M}_{0,n+\Sigma_{i}n_{i}}$ is the moduli space of 
genus zero stable $n+\Sigma_{i}n_{i}$-marked curves, then the $\psi$ classes on $\overline{\mathcal{M}}_{n+\Sigma_{i}n_{i}}(\mathcal{X},0,\mathbf{x})$
is defined to be the pullback of the $\psi$ classes on $\overline{M}_{0,n+\Sigma_{i}n_{i}}$. As in Section \ref{review}, the degree zero and genus zero descendent orbifold Gromov-Witten invariants are defined by
\begin{equation}\Label{orbinv}
\langle\tau_{l_{1}}\cdots\tau_{l_{n}}\widetilde{\tau}_{k_{1}}\cdots\widetilde{\tau}_{k_{\Sigma_{i}
n_{i}}}\rangle^{\mathcal{X}}
:=\int_{[\overline{\mathcal{M}}_{n+\Sigma_{i}n_{i}}(\mathcal{X},0,\mathbf{x})]^{vir}}
\psi^{l_{1}}\cdots\psi^{l_{n}}\psi^{k_{1}}\cdots\psi^{k_{\Sigma_{i}
n_{i}}}.
\end{equation}

Let $L_\omega$ denote the line bundle over $\mathcal{B}\zz_N$ defined by the $\zz_N$-representation $\cc$ on which $\omega\in \zz_N$ acts by multiplication by $\omega$. The twisted stable maps factor through $\mathcal{B}\mathbb{Z}_{N}$ and $\mathcal{B}\mathbb{Z}_{N}\subset \mathcal{X}$ has normal bundle $L_{\omega}^{\oplus d}$ of rank $d$ since the action of $\mathbb{Z}_{N}$ on the neighborhood of $p$ is given by the diagonal matrix $diag(\omega,\cdots,\omega)$. Consider the following diagram:
$$\xymatrix{
~\widetilde{\mathcal{C}}\dto_{p}\rto & pt\dto^{} \\
~\mathcal{C}\rto^{f}\dto_{\pi} &  ~\mathcal{B}\mathbb{Z}_{N}\\
\overline{\mathcal{M}}_{n+\Sigma_{i}
n_{i}}(\mathcal{B}\mathbb{Z}_{N}, \mathbf{x}).&~}$$ 

Let $$\widetilde{\pi}=\pi\circ p: \widetilde{\mathcal{C}}\to \overline{\mathcal{M}}_{n+\Sigma_{i}n_{i}}(\mathcal{B}\mathbb{Z}_{N},\mathbf{x})$$
be the composite map. Let $\mathbb{E}^{\vee}=R^{1}\widetilde{\pi}_{*}\mathcal{O}$ be the dual Hodge bundle. The action of $\omega\in\mathbb{Z}_{N}$ on $\widetilde{\mathcal{C}}$ induces an action of $\omega$ on $\mathbb{E}^{\vee}$, giving the decomposition into eigen-bundles:
$$\mathbb{E}^{\vee}=\mathbb{E}_{1}^{\vee}\oplus \mathbb{E}_{\omega}^{\vee}
\oplus\cdots\oplus \mathbb{E}_{\omega^{N-1}}^{\vee}.$$
Following the convention in \cite{BGP}, $\mathbb{E}_{\omega^i}^\vee$ is the eigen-bundle on which $\omega$ acts with eigenvalue $\omega^i$.

It is easy to see that
$$\mathbb{E}_{1}^{\vee}=0,~~R^{1}\pi_{*}f^{*}(L_{\omega})=\mathbb{E}_{\omega^{N-1}}^{\vee}.$$
So
\begin{align}\Label{orbinv2}
\langle\tau_{l_{1}}\cdots\tau_{l_{n}}\widetilde{\tau}_{k_{1}}\cdots\widetilde{\tau}_{k_{\Sigma_{i}
n_{i}}}\rangle^{\mathcal{X}}
&=\int_{\overline{\mathcal{M}}_{n+\Sigma_{i}n_{i}}(\mathcal{B}\mathbb{Z}_{N}, \mathbf{x})}
\psi^{l_{1}}\cdots\psi^{l_{n}}\psi^{k_{1}}\cdots\psi^{k_{\Sigma_{i}
n_{i}}} e(\mathbb{E}_{\omega^{N-1}}^{\vee}\oplus\cdots\oplus
\mathbb{E}_{\omega^{N-1}}^{\vee})\nonumber \\
&=\int_{\overline{\mathcal{M}}_{n+\Sigma_{i}n_{i}}(\mathcal{B}\mathbb{Z}_{N}, \mathbf{x})}
\psi^{l_{1}}\cdots\psi^{l_{n}}\psi^{k_{1}}\cdots\psi^{k_{\Sigma_{i}
n_{i}}} \lambda_{r_{1}}^{d},
\end{align}
where $e$ is the  Euler class and $\lambda_{r_{1}}$ the top Chern class of the bundle
$\mathbb{E}_{\omega}\cong\mathbb{E}_{\omega^{N-1}}^{\vee}$.  Let
$r_{1}=rank(\mathbb{E}_{\omega})$,
$r_{2}=rank(\mathbb{E}_{\omega^{2}}),\cdots$,
$r_{N-1}=rank(\mathbb{E}_{\omega^{N-1}})$, then
\begin{prop}\Label{dimension} We have
$$r_{1}=\sum_{i=1}^{N-1}n_i\frac{i}{N} -1, \quad r_{N-1}=\sum_{i=1}^{N-1}n_i\frac{N-i}{N}-1,\quad 
\text{and } r_{1}+r_{N-1}-1=\Sigma_{i}n_{i}-3.$$
\end{prop}
\begin{proof}
Let $f: \C\to \B\zz_N$ be an orbifold stable map with stack structures described by (\ref{ramification}). By Riemann-Roch for stacky curves, we find 
\begin{equation*}
\chi(\C, f^*L_\omega)=1-\sum_{i=1}^{N-1}n_i\frac{i}{N}, \quad \chi(\C, f^*L_{\omega^{N-1}})=1-\sum_{i=1}^{N-1}n_i\frac{N-i}{N}.
\end{equation*}
We conclude by observing that 
\begin{equation*}
\begin{split}
&R^0\pi_*f^*L_\omega=0, \quad R^1\pi_*f^*L_\omega=\mathbb{E}_\omega;\\
&R^0\pi_*f^*L_{\omega^{N-1}}=0, \quad R^1\pi_*f^*L_{\omega^{N-1}}=\mathbb{E}_{\omega^{N-1}}.
\end{split}
\end{equation*}
\end{proof}

In the next two sections we consider the case for the weighted projective line $\mathbb{P}(1,N)$ and the weighted projective surface $\mathbb{P}(1,1,N)$. As discussed above, every stable map in
$\overline{\mathcal{M}}_{n+\Sigma_{i}n_{i}}(\mathcal{X},0,\mathbf{x})$ determines (and is determined by) a $\mathbb{Z}_{N}$-admissible cover over the coarse curve $\widetilde{\mathcal{C}}\to C$, with ramification data specified by $\mathbf{x}$ in (\ref{ramification}). The genus, denoted by $g$, of the cover $\widetilde{\mathcal{C}}$ is fixed by $\mathbf{x}$. Introduce
the following notations for the descendent degree zero genus zero orbifold Gromov-Witten
invariants of $\mathbb{P}(1,N)$ and $\mathbb{P}(1,1,N)$:

\begin{equation}\Label{orbinv3}
\langle\tau_{l_{1}}\cdots\tau_{l_{n}}\widetilde{\tau}_{k_{1}}\cdots\widetilde{\tau}_{k_{\Sigma_{i}
n_{i}}}|\lambda_{r_{1}}\rangle_{g}
=\int_{\overline{\mathcal{M}}_{n+\Sigma_{i}n_{i}}(\mathcal{B}\mathbb{Z}_{N},\mathbf{x})}
\psi^{l_{1}}\cdots\psi^{l_{n}}\psi^{k_{1}}\cdots\psi^{k_{\Sigma_{i}
n_{i}}} \lambda_{r_{1}},
\end{equation}

\begin{equation}\Label{orbinv4}
\langle\tau_{l_{1}}\cdots\tau_{l_{n}}\widetilde{\tau}_{k_{1}}\cdots\widetilde{\tau}_{k_{\Sigma_{i}
n_{i}}}|\lambda^{2}_{r_{1}}\rangle_{g}
=\int_{\overline{\mathcal{M}}_{n+\Sigma_{i}n_{i}}(\mathcal{B}\mathbb{Z}_{N}, \mathbf{x})}
\psi^{l_{1}}\cdots\psi^{l_{n}}\psi^{k_{1}}\cdots\psi^{k_{\Sigma_{i}
n_{i}}} \lambda^{2}_{r_{1}}
\end{equation}
and
\begin{multline}\label{orbinv5}
\langle\langle\tau_{l_{1}}\cdots\tau_{l_{n}}\widetilde{\tau}_{k_{1}}\cdots\widetilde{\tau}_{k_{\Sigma_{i}
n_{i}}}|\lambda_{r_{1}}\rangle\rangle_{g}=\\ \sum_{M\geq
0}\frac{1}{M!}\sum_{b_1,\cdots,b_M\geq 0}t_{b_{1}}\cdots
t_{b_{M}}\langle\tau_{b_{1}}\cdots\tau_{b_{M}}\tau_{l_{1}}\cdots\tau_{l_{n}}\widetilde{\tau}_{k_{1}}\cdots\widetilde{\tau}_{k_{\Sigma_{i}
n_{i}}}|\lambda_{r_{1}}\rangle_{g},
\end{multline}

\begin{multline}\label{orbinv6}
\langle\langle\tau_{l_{1}}\cdots\tau_{l_{n}}\widetilde{\tau}_{k_{1}}\cdots\widetilde{\tau}_{k_{\Sigma_{i}
n_{i}}}|\lambda^{2}_{r_{1}}\rangle\rangle_{g}=\\ \sum_{M\geq
0}\frac{1}{M!}\sum_{b_1,\cdots,b_M\geq 0}t_{b_{1}}\cdots
t_{b_{M}}\langle\tau_{b_{1}}\cdots\tau_{b_{M}}\tau_{l_{1}}\cdots\tau_{l_{n}}\widetilde{\tau}_{k_{1}}\cdots\widetilde{\tau}_{k_{\Sigma_{i}
n_{i}}}|\lambda^{2}_{r_{1}}\rangle_{g}.
\end{multline}

\begin{rmk}
If $n=0$, then the moduli stack
$\overline{\mathcal{M}}_{\Sigma_{i}n_{i}}(\mathcal{B}\mathbb{Z}_{N}, \mathbf{x})$
is the Hurwitz scheme $\overline{H}_{g}$ parametrizing the
admissible $\mathbb{Z}_{N}$-covers to $\mathbb{P}^{1}$ with
ramification type $\mathbf{x}$ in (\ref{ramification}). 
The
integrals above are called Hurwitz-Hodge integrals in \cite{BGP} (see also \cite{BJ}). In
the following we always use this monodromy data for an admissible
$\mathbb{Z}_{N}$-covers to $\mathbb{P}^{1}$.
\end{rmk}

\section{The degree zero  Virasoro conjecture for weighted projective lines.}\label{weighted_P1}

Let $\mathcal{X}=\mathbb{P}(1,N)$ be the weighted projective space
with weights $1, N$ for $N\in \mathbb{Z}_{>0}$. The Chen-Ruan orbifold
cohomology $H^{*}_{CR}(\mathbb{P}(1,N))$ has the following
generators: $$1\in H^{0}_{CR}(\mathbb{P}(1,N)), \,
\xi \in H^{2}_{CR}(\mathbb{P}(1,N)), \,
\gamma_{j}\in H^{\frac{2j}{N}}_{CR}(\mathbb{P}(1,N)), \,
\text{ for }1\leq j\leq N-1.$$

We will use the following coordinates for the corresponding descendents:
\begin{center}
\begin{tabular}{lp{0.75\textwidth}}
coordinate & \quad descendent\\
$t_i, \, i\geq 0$ & \quad $\tau_i(1)$\\
$s_{i}, \, i\geq 0$ & \quad $\tau_{i}(\xi)$\\
$\alpha^{j}_{i}, \, i\geq 0$ & \quad $\tau_{i}(\gamma_{j})$ for $1\leq j\leq N-1$.
\end{tabular}
\end{center}

Let $\mathcal{F}_{\mathcal{X}}^{0,0}$ be the degree zero, genus zero
orbifold Gromov-Witten potential of $\mathcal{X}$. Let 
$$\mathcal{D}_{\mathcal{X}}^{0,0}=\exp\left(\frac{1}{\hbar}\mathcal{F}_{\mathcal{X}}^{0,0}\right).$$ 

For $1\leq i\leq N-1$ let 
\begin{equation}\label{notation}
M_{i}:=\sum_{a=0}^{i}n_{a},
\end{equation} 
where $n_{0}:=0$. This notation is used in Section \ref{weighted_P1} and \ref{weighted_P2}.

\begin{thm}
We have
$$\frac{\mathcal{L}_{k}\mathcal{D}_{\mathcal{X}}^{0,0}}{\mathcal{D}_{\mathcal{X}}^{0,0}}=
\frac{1}{\hbar}\left((-2)x_{0,l}^{k}(\mathbf{t})+\sum_{l=0}^{\infty}s_{l}\cdot
y_{0,l}^{k}(\mathbf{t})
+\sum_{g=0}^{\infty}~\sum_{k_{1},\cdots,k_{\Sigma
n_{i}}\geq 0}\alpha_{k_{1}}^{1}\cdots\alpha_{k_{\Sigma_{i}n_{i}}}^{N-1}\cdot
z_{g;k_{1},\cdots,k_{\Sigma_{i}n_{i}}}^{k}(\mathbf{t})\right),$$ where

\begin{multline*}
x_{0,l}^{k}(\mathbf{t})=
-[1]_{0}^{k}\langle\langle\tau_{k+1}|\lambda_{0}\rangle\rangle_{0}
+\sum_{m=1}^{\infty}[m]_{0}^{k}t_{m}\langle\langle\tau_{k+m}|\lambda_{0}\rangle\rangle_{0}
+ [1]_{1}^{k}\langle\langle\tau_{k}|\lambda_{0}\rangle\rangle_{0} \\
-\sum_{m=0}^{\infty}[m]_{1}^{k}t_{m}\langle\langle\tau_{k+m-1}|\lambda_{0}\rangle\rangle_{0}
-\frac{1}{2}\sum_{m=0}^{k-2}(-1)^{m+1}[-m-1]_{1}^{k}\langle\langle\tau_{m}|\lambda_{0}\rangle\rangle_{0}
\langle\langle\tau_{k-m-2}|\lambda_{0}\rangle\rangle_{0},
\end{multline*}
\begin{equation*}
y_{0,l}^{k}(\mathbf{t})=
-[1]_{0}^{k}\langle\langle\tau_{k+1}\tau_{l}|\lambda_{0}\rangle\rangle_{0}
+\sum_{m=1}^{\infty}[m]_{0}^{k}t_{m}\langle\langle\tau_{k+m}\tau_{l}|\lambda_{0}\rangle\rangle_{0}
+
[l+1]_{0}^{k}\langle\langle\tau_{k+l}|\lambda_{0}\rangle\rangle_{0}
\end{equation*}
and
\begin{multline*}
z_{g;k_{1},\cdots,k_{\Sigma_{i}n_{i}}}^{k}(\mathbf{t})=
-[1]_{0}^{k}\langle\langle\tau_{k+1}\widetilde{\tau}_{k_{1}}\cdots
\widetilde{\tau}_{k_{\Sigma_{i}n_{i}}}|\lambda_{r_{1}}\rangle\rangle_{g}
+\sum_{m=1}^{\infty}[m]_{0}^{k}t_{m}\langle\langle\tau_{k+m}\widetilde{\tau}_{k_{1}}\cdots \widetilde{\tau}_{k_{\Sigma_{i}n_{i}}}|\lambda_{r_{1}}\rangle\rangle_{g}\\
+ \sum_{i=1}^{N-1}\sum_{j=M_{i-1}+1}^{M_{i}}[k_{j}+\tfrac{i}{N}]_{0}^{k}
\langle\langle\widetilde{\tau}_{k_{1}}\cdots\widetilde{\tau}_{k_{n_{i-1}}}\widetilde{\tau}_{k_{n_{i-1}+1}}\cdots\widetilde{\tau}_{k+k_{j}}\cdots\widetilde{\tau}_{k_{n_{i}}}\widetilde{\tau}_{k_{n_{i}+1}}
\cdots\widetilde{\tau}_{k_{\Sigma_{i}n_{i}}}|\lambda_{r_{1}}\rangle\rangle_{g}.
\end{multline*}
\end{thm}

\begin{proof}
The Virasoro operator $\mathcal{L}_{k}$, $k>0$ is given by
\begin{multline*}
\mathcal{L}_{k}= -[1]_{0}^{k}\partial
t_{k+1}+\sum_{m=0}^{\infty}\Bigl([m]_{0}^{k}t_{m}\partial t_{k+m}
+[m+1]_{0}^{k}s_{m}\partial s_{k+m}+\sum_{i=1}^{N-1}[m+\tfrac{i}{N}]_{0}^{k}\alpha_{m}^{i}\partial \alpha_{k+m}^{i} 
\Bigr)\\ +
2\Bigl(-[1]_{1}^{k}\partial s_{k}+\sum_{m=0}^{\infty}[m]_{1}^{k}t_{m}\partial s_{k+m-1} 
+\frac{\hbar}{2}\sum_{m=0}^{k-2}(-1)^{m+1}[-m-1]_{1}^{k}\partial
s_{m}\partial s_{k-m-2}\Bigr).
\end{multline*}
The degree zero orbifold Gromov-Witten potential is
\begin{multline*}
\mathcal{D}_{\mathcal{X}}^{0,0}=
\exp\Bigl(\frac{1}{\hbar}\sum_{m=0}^{\infty}s_{m}\langle\langle\tau_{m}|\lambda_{0}\rangle\rangle_{0}
+\frac{1}{\hbar}(-2)\langle\langle 1\rangle\rangle_{0} \\
+\frac{1}{\hbar}\sum_{g=0}^{\infty}\sum_{k_{1},\cdots,k_{\Sigma_{i}n_{i}}\geq 0}\alpha^{1}_{k_{1}}\cdots
\alpha^{N-1}_{k_{\Sigma_{i}n_{i}}}\langle\langle\widetilde{\tau}_{k_{1}}\cdots
\widetilde{\tau}_{k_{\Sigma_{i}n_{i}}}|\lambda_{r_{1}}\rangle\rangle_{g}\Bigr).
\end{multline*}
Applying the operator $\mathcal{L}_{k}$ to
$\mathcal{D}_{\mathcal{X}}^{0,0}$ we obtain the result.
\end{proof}

So the degree zero, genus zero Virasoro constraints for $\mathbb{P}(1,N)$ is equivalent to the vanishing of $x_{0,l}^{k}(\mathbf{t}), y_{0,l}^{k}(\mathbf{t})$ and $z_{g,k_{1},\cdots,k_{\Sigma_{i}n_{i}}}^{k}(\mathbf{t})$.

Now let
\begin{equation}\Label{constant1}
\Gamma_{j,g}=\langle\widetilde{\tau}_{a}|\lambda_{r_{1}}\rangle_{g}=\int_{\overline{H}_{g}}\psi^{a}\lambda_{r_{1}}, \quad \text{where } a:=\sum_{i=1}^{N-1}n_{i}-\sum_{i=1}^{N-1}\frac{i}{N}n_{i}-2.
\end{equation}

Write $\mathbf{\Gamma_{g}}:=(\Gamma_{j,g})_{1\leq j\leq \Sigma_{i}n_{i}}$ as a column vector. Let $\mathbf{c_{g}}:=(c_{j,g})_{1\leq j\leq \Sigma_{i}n_{i}}$ be another column vector. Let the index $i$ vary from $1$ to $N-1$ and define a $\Sigma_{i}n_{i}\times\Sigma_{i}n_{i}$ square matrix $A=(a_{st})$ by

\begin{equation}\Label{matrix}
a_{st}:=\begin{cases}\frac{i}{N}+a&\text{if
}~ M_{i-1}<s=t\leq M_{i}\,;\\
\frac{i}{N}&\text{if}~M_{i-1}<t\leq M_{i}~\text{and}~s\neq t\,.\end{cases}
\end{equation}

The matrix can be written as follows:
$$
A=\left[
\begin{array}{cccccccccc}
  \tfrac{1}{N}+a&\cdots&\tfrac{1}{N}&\tfrac{2}{N}&\cdots&\tfrac{2}{N}&\cdots&\tfrac{N-1}{N}&\cdots&\tfrac{N-1}{N} \\
  \tfrac{1}{N}&\cdots&\tfrac{1}{N}&\tfrac{2}{N}&\cdots&\tfrac{2}{N}&\cdots&\tfrac{N-1}{N}&\cdots&\tfrac{N-1}{N} \\
  \cdots&\cdots&\cdots&\cdots&\cdots&\cdots&\cdots&\cdots&\cdots&\cdots\\
  \tfrac{1}{N}&\cdots&\tfrac{1}{N}+a&\tfrac{2}{N}&\cdots&\tfrac{2}{N}&\cdots&\tfrac{N-1}{N}&\cdots&\tfrac{N-1}{N} \\
  \tfrac{1}{N}&\cdots&\tfrac{1}{N}&\tfrac{2}{N}+a&\cdots&\tfrac{2}{N}&\cdots&\tfrac{N-1}{N}&\cdots&\tfrac{N-1}{N} \\
  \cdots&\cdots&\cdots&\cdots&\cdots&\cdots&\cdots&\cdots&\cdots&\cdots \\
  \tfrac{1}{N}&\cdots&\tfrac{1}{N}&\tfrac{2}{N}&\cdots&\tfrac{2}{N}+a&\cdots&\tfrac{N-1}{N}&\cdots&\tfrac{N-1}{N} \\
 \cdots&\cdots&\cdots&\cdots&\cdots&\cdots&\cdots&\cdots&\cdots&\cdots \\
 \tfrac{1}{N}&\cdots&\tfrac{1}{N}&\tfrac{2}{N}&\cdots&\tfrac{2}{N}&\cdots&\tfrac{N-1}{N}+a&\cdots&\tfrac{N-1}{N}\\
  \cdots&\cdots&\cdots&\cdots&\cdots&\cdots&\cdots&\cdots&\cdots&\cdots \\
 \tfrac{1}{N}&\cdots&\tfrac{1}{N}&\tfrac{2}{N}&\cdots&\tfrac{2}{N}&\cdots&\tfrac{N-1}{N}&\cdots&\tfrac{N-1}{N}+a
\end{array}
\right].
$$
It is easy to check that $A$ is nonsingular for $a\neq 0$. Let
$\mathbf{A}$ be the matrix obtained from $A$ as follows: for an integer $j$ with $M_{i-1}+1\leq j\leq M_i$ for some $1\leq i\leq N-1$, the $j$-th row of $\mathbf{A}$ is obtained by multiplying the $j$-th row of $A$ by $$\frac{(\frac{i}{N})((\sum_{i=1}^{N-1}n_i)-3)!}{(a+\frac{i}{N})!\prod_{i=1}^{N-1}(\frac{i}{N})^{n_i}}.$$
Here $(a+\frac{i}{N})!:=\prod_{m=0}^a(m+\frac{i}{N}).$


The linear system
\begin{equation}\label{linearsystem}
\mathbf{A}\cdot \mathbf{c_{g}}=\mathbf{\Gamma_{g}}
\end{equation}
has a unique solution which represents $c_{j,g}$ as a linear combination of $\Gamma_{j,g}$'s for $1\leq j\leq \Sigma_{i}n_{i}$.

For integers $1\leq s\leq N-1$ and $1\leq j\leq\Sigma_{i}n_{i}$, we put 
$$(k_{j}+\frac{s}{N})!=\frac{s}{N}\cdot(1+\frac{s}{N})\cdots(k_{j}+\frac{s}{N}).$$

From the vanishing of $z_{g,k_{1},\cdots,k_{\Sigma_{i}n_{i}}}^{k}(\mathbf{t})$ for $k\geq 1$ and $k_{1},\cdots,k_{\Sigma_{i}n_{i}}\geq 0$, we obtain the following theorem.

\begin{thm}\Label{curvevir}
We have
$$
\langle\widetilde{\tau}_{k_{1}}\cdots\widetilde{\tau}_{k_{\Sigma_{i}n_{i}}}\tau_{l_{1}}
\cdots\tau_{l_{n}}|\lambda_{r_{1}}\rangle_{g}=
\sum_{s=1}^{N-1}\sum_{j=M_{s-1}+1}^{M_{s}}\tfrac{(n+\Sigma_{i}
n_{i}-3)!(k_{j}+\frac{s}{N})}
{\prod_{j}l_{j}!\prod_{b=1}^{N-1}\prod_{j=M_{b-1}+1}^{M_{b}}(k_{j}+\frac{b}{N})!}c_{j,g}.
$$
\end{thm}

\begin{proof}
The vanishing of
$z_{g,k_{1},\cdots,k_{\Sigma_{i}n_{i}}}^{k}(\mathbf{t})$ gives the vanishing of its Taylor coefficients, $$\tfrac{1}{[1]_{0}^{k}}\partial t_{l_{1}}\cdots \partial
t_{l_{n}}z_{g;k_{1},\cdots,k_{\Sigma_{i}n_{i}}}^{k}(0)=0.$$
An explicit calculation shows that $\tfrac{1}{[1]_{0}^{k}}\partial t_{l_{1}}\cdots \partial
t_{l_{n}}z_{g;k_{1},\cdots,k_{\Sigma_{i}n_{i}}}^{k}(0)$ is the right side of the following
\begin{multline}\Label{recursion}
0=-\langle\tau_{k+1}\widetilde{\tau}_{k_{1}}\cdots\widetilde{\tau}_{k_{\Sigma_{i}n_{i}}}\tau_{l_{1}}
\cdots\tau_{l_{n}}|\lambda_{r_{1}}\rangle_{g}+\sum_{i=1}^{n}\tfrac{(l_{i}+k)!}{(l_{i}-1)!(k+1)!}\langle\widetilde{\tau}_{k_{1}}\cdots\widetilde{\tau}_{k_{\Sigma_{i}n_{i}}}\tau_{l_{1}}
\cdots\tau_{l_{i}+k}\cdots\tau_{l_{n}}|\lambda_{r_{1}}\rangle_{g}\\
+
\sum_{i=1}^{N-1}\sum_{j=M_{i-1}+1}^{M_{i}}\tfrac{(k_{j}+k+\frac{i}{N})!}{(k_{j}-1+\frac{i}{N})!(k+1)!}
\langle\widetilde{\tau}_{k_{1}}\cdots\widetilde{\tau}_{k_{n_{i-1}}}\widetilde{\tau}_{k_{n_{i-1}+1}}\cdots\widetilde{\tau}_{k_{j}+k}\cdots
\widetilde{\tau}_{k_{n_{i}}}\widetilde{\tau}_{k_{n_{i}+1}}\cdots\widetilde{\tau}_{k_{\Sigma_{i}n_{i}}}\tau_{l_{1}}
\cdots\tau_{l_{n}}|\lambda_{r_{1}}\rangle_{g}.
\end{multline}

We now solve the recursion (\ref{recursion}). The virtual dimension of $\overline{\mathcal{M}}_{n+1+\Sigma_{i}n_{i}}(\mathbb{P}(1,N),0,\mathbf{x})$ is

$$\text{vdim}=1+n+\sum_{i=1}^{N-1}n_{i}+1-3-\sum_{i=1}^{N-1}n_{i}\frac{i}{N}.$$

If $\langle\tau_{k+1}\widetilde{\tau}_{k_{1}}\cdots\widetilde{\tau}_{k_{\Sigma_{i}n_{i}}}\tau_{l_{1}}
\cdots\tau_{l_{n}}|\lambda_{r_{1}}\rangle_{g}\neq 0$, we have
$\text{vdim}=k+1+\sum_{i=1}^{n}l_{i}+\sum_{i=1}^{\Sigma_{i}n_{i}}k_{i}$. So

\begin{equation}\Label{virdim}
n+\sum_{i=1}^{N-1}
n_{i}-2=\sum_{i=1}^{n}l_{i}+\sum_{i=1}^{N-1}\sum_{j=M_{i-1}+1}^{M_{i}}\Bigl(k_{j}+\tfrac{i}{N}\Bigr)+k.
\end{equation}

For an integer $r$, then there exists a unique  integer $s$ such that $M_{s-1}+1\leq r\leq M_{s}$, with $1\leq s\leq N-1$. 
Let $\mathbf{k}=(k_{1},\cdots,k_{\Sigma_{i}n_{i}})$ and 
$\mathbf{l}=(l_{1},\cdots,l_{n})$. Introduce 

\begin{equation}\label{special}
\Theta(\mathbf{k},\mathbf{l})_{r}:=\frac{(n+\Sigma_{i}
n_{i}-3)!(k_{r}+\frac{s}{N})}
{\prod_{j}l_{j}!\prod_{b=1}^{N-1}\prod_{j=M_{b-1}+1}^{M_{b}}(k_{j}+\frac{b}{N})!}.
\end{equation}

We claim that $\Theta(\mathbf{k},\mathbf{l})_{r}$ is a solution of the recursion (\ref{recursion}). To see this, write (\ref{virdim})
as
\begin{equation}\Label{virdimm}
n+\sum_{i=1}^{N-1}
n_{i}-2=\sum_{i=1}^{n}l_{i}+\sum_{i=1}^{N-1}\sum_{j=M_{i-1}+1,j\neq r}^{M_{i}}\Bigl(k_{j}+\tfrac{i}{N}\Bigr)+(k_{r}+\frac{s}{N}+k).
\end{equation}
Multiply on both sides of (\ref{virdimm}) by
$$
\frac{(n+\Sigma_{i}n_{i}-3)!(k_{r}+\tfrac{s}{N})}
{(k+1)!\prod_{j}l_{j}!\prod_{b=1}^{N-1}\prod_{j=M_{b-1}+1}^{M_{b}}(k_{j}+\frac{b}{N})!},
$$
we get
\begin{equation*}
\begin{split}
&\frac{(n+\Sigma_{i}n_{i}-2)!(k_{r}+\frac{s}{N})}
{(k+1)!\prod_{j}l_{j}!\prod_{b=1}^{N-1}\prod_{j=M_{b-1}+1}^{M_{b}}(k_{j}+\frac{b}{N})!}\\
&= \sum_{i=1}^{n}\frac{(l_{i}+k)!}{(l_{i}-1)!(k+1)!}\frac{(n+\Sigma_{i}n_{i}-3)!(k_{r}+\frac{s}{N})}
{l_{1}!\cdots (l_{i}+k)!\cdots l_{n}!\prod_{b=1}^{N-1}\prod_{j=M_{b-1}+1}^{M_{b}}(k_{j}+\frac{b}{N})!} \\
&+\sum_{i=1}^{N-1}\frac{(k_{j}+k+\frac{i}{N})!}{(k_{j}-1+\frac{i}{N})!(k+1)!}\cdot
\sum_{j=M_{i-1}+1,j\neq r}^{M_{i}} \tfrac{(n+\Sigma_{i}n_{i}-3)!(k_{r}+\frac{s}{N})}
{\prod_{j}l_{j}!\prod_{b=1,b\neq i}^{N-1}\prod_{j=M_{b-1}+1}^{M_{b}}(k_{j}+\frac{b}{N})!(k_{n_{i-1}+1}+\frac{i}{N})!\cdots (k_{j}+k+\frac{i}{N})!\cdots(k_{n_{i}}+\frac{i}{N})!}\\
&+\frac{(k_{r}+k+\frac{s}{N})!}{(k_{r}-1+\frac{s}{N})!(k+1)!}\cdot \tfrac{(n+\Sigma_{i}n_{i}-3)!(k_{r}+k+\frac{s}{N})}{\prod_{j}l_{j}!\prod_{b=1,b\neq s}^{N-1}\prod_{j=M_{b-1}+1}^{M_{b}}(k_{j}+\frac{b}{N})!(k_{n_{s-1}+1}+\frac{s}{N})!\cdots (k_{r}+k+\frac{s}{N})!\cdots(k_{n_{s}}+\frac{i}{N})!}.
\end{split}
\end{equation*}
It is straightforward to see that this is the recursion (\ref{recursion}). 

Suppose that $\langle\widetilde{\tau}_{k_{1}}\cdots\widetilde{\tau}_{k_{\Sigma_{i}n_{i}}}\tau_{l_{1}}
\cdots\tau_{l_{n}}|\lambda_{r_{1}}\rangle_{g}$ is of the form $\sum_r c_{r,g}\Theta(\mathbf{k},\mathbf{l})_r$. Then by considering special values of $\mathbf{k}, \mathbf{l}$, we find that the coefficients $c_{r,g}$ are uniquely determined by the linear system (\ref{linearsystem}). The result follows.

\end{proof}

The initial values in Theorem \ref{curvevir} are the following Hurwitz-Hodge integrals
$$\int_{\overline{H}_{g}}\lambda_{r_{1}}\psi^{r_{N-1}-1},$$
where $\overline{H}_{g}$ is the Hurwitz scheme parametrizing admissible $\mathbb{Z}_{N}$-covers over $\mathbb{P}^{1}$. These integrals can be interpreted as certain orbifold Gromov-Witten invariants of $\mathcal{B}\mathbb{Z}_N$ twisted by the line bundle $L_\omega$ and the inverse (equivariant) Euler class. The main results of \cite{ccit2} can be applied to compute these integrals.

\section{The degree zero Virasoro conjecture for weighted projective surfaces.}\label{weighted_P2}

Let $\X=\mathbb{P}(1,1,N)$ be the weighted projective plane with weights $1,1, N$ for
$N\in \mathbb{Z}_{>0}$. The Chen-Ruan orbifold cohomology $H^{*}_{CR}(\mathbb{P}(1,1,N))$
has generators as follows:
\begin{equation*}
\begin{split}
& 1\in H^{0}_{CR}(\mathbb{P}(1,1,N)), \,\xi\in H^{2}_{CR}(\mathbb{P}(1,1,N), \, [\mathcal{X}]\in H^{4}_{CR}(\mathbb{P}(1,1,N)),\\
& \gamma_{j}\in H^{\tfrac{4j}{N}}_{CR}(\mathbb{P}(1,1,N)) \, \text{ for }1\leq j\leq N-1.
\end{split}
\end{equation*}

We will use the following coordinates for the corresponding descendants:
\begin{center}
\begin{tabular}{lp{0.75\textwidth}}
coordinate & \quad descendent\\
$t_i, \, i\geq 0$ & \quad $\tau_i(1)$\\
$s_{i}, \, i\geq 0$ & \quad $\tau_{i}(\xi)$\\
$r_{i}, \, i\geq 0$ & \quad $\tau_{i}([\mathcal{X}])$\\
$\alpha^{j}_{i}, \, i\geq 0$ & \quad $\tau_{i}(\gamma_{j})$ for $1\leq j\leq N-1$.
\end{tabular}
\end{center}

Let $\mathcal{F}_{\mathcal{X}}^{0,0}$ be the degree zero, genus zero
orbifold Gromov-Witten potential of $\mathcal{X}$. Let
$$\mathcal{D}_{\mathcal{X}}^{0,0}=\exp\left(\frac{1}{\hbar}\mathcal{F}_{\mathcal{X}}^{0,0}\right).$$

\begin{thm}
We have
$$\frac{\mathcal{L}_{k}\mathcal{D}_{\mathcal{X}}^{0,0}}{\mathcal{D}_{\mathcal{X}}^{0,0}}=
\frac{1}{\hbar}w(r,s,t)
+\frac{1}{\hbar}\sum_{g=0}^{\infty}~\sum_{k_{1},\cdots,k_{\Sigma_{i}n_{i}}\geq 0}\alpha^{1}_{k_{1}}\cdots\alpha^{N-1}_{k_{\Sigma_{i}
n_{i}}}\cdot
z_{g;k_{1},\cdots,k_{\Sigma_{i}n_{i}}}^{k}(\mathbf{t}),$$ where
$w(r,s,t)$ is  the general degree zero, genus zero potential without
stacky points on the curve in \cite{GP}, and
\begin{multline*}
z_{g;k_{1},\cdots,k_{\Sigma_{i}n_{i}}}^{k}(\mathbf{t})=
-[\tfrac{1}{2}]_{0}^{k}\langle\langle\tau_{k+1}\widetilde{\tau}_{k_{1}}\cdots
\widetilde{\tau}_{k_{\Sigma_{i}n_{i}}}|\lambda^{2}_{r_{1}}\rangle\rangle_{g}
+\sum_{m=1}^{\infty}[m-\tfrac{1}{2}]_{0}^{k}t_{m}\langle\langle\tau_{k+m}\widetilde{\tau}_{k_{1}}\cdots \widetilde{\tau}_{k_{\Sigma_{i}n_{i}}}|\lambda^{2}_{r_{1}}\rangle\rangle_{g}\\
+\sum_{i=1}^{N-1}\sum_{j=M_{i-1}+1}^{M_{i}}[k_{j}+\tfrac{2i}{N}-\tfrac{1}{2}]_{0}^{k}
\langle\langle\widetilde{\tau}_{k_{1}}\cdots\widetilde{\tau}_{k_{n_{i-1}}}\widetilde{\tau}_{k_{n_{i-1}+1}}\cdots\widetilde{\tau}_{k+k_{j}}\cdots\widetilde{\tau}_{k_{n_{i}}}\widetilde{\tau}_{k_{n_{i}+1}}
\cdots\widetilde{\tau}_{k_{\Sigma_{i}n_{i}}}|\lambda^{2}_{r_{1}}\rangle\rangle_{g}.
\end{multline*}
\end{thm}

\begin{proof}
The Virasoro operator $\mathcal{L}_{k}$, $k>0$ is given by
\begin{align*}
\mathcal{L}_{k}&= -[\tfrac{1}{2}]_{0}^{k}\partial
t_{k+1}+\sum_{m=0}^{\infty}\Bigl([m-\tfrac{1}{2}]_{0}^{k}t_{m}\partial
t_{k+m}
+[m+\tfrac{1}{2}]_{0}^{k}s_{m}\partial s_{k+m}+[m+\tfrac{3}{2}]_{0}^{k}r_{m}\partial r_{k+m} \\
&+\sum_{i=1}^{N-1}[m+\tfrac{2i}{N}-\tfrac{1}{2}]_{0}^{k}\alpha_{m}^{i}\partial
\alpha_{k+m}^{i}\Bigr)
+ \hbar\sum_{m=0}^{k-1}(-1)^{m+1}\Bigl([-m-\tfrac{3}{2}]_{1}^{k}\partial
r_{m}\partial t_{k-m-1}+\\
&\tfrac{1}{2}[-m-\tfrac{1}{2}]_{1}^{k}\partial s_{m}\partial s_{k-m-1}
+ \sum_{i=1}^{N-1}[-m-\tfrac{2i}{N}+\tfrac{1}{2}]_{1}^{k}\partial
\alpha^{i}_{m}\partial \alpha^{N-i}_{k-m-1}\Bigr) + \\
&\mathbf{c}\cdot\Bigl(-[\tfrac{1}{2}]_{1}^{k}\partial s_{k}+\sum_{m=0}^{\infty}[m-\tfrac{1}{2}]_{1}^{k}t_{m}\partial s_{k+m-1} 
+[m+\tfrac{1}{2}]_{1}^{k}s_{m}\partial
r_{k+m-1}\\
&+\hbar\sum_{m=0}^{k-2}(-1)^{m+1}[-m-\tfrac{3}{2}]_{1}^{k}\partial
r_{m}\partial s_{k-m-2}
\Bigr)
+|\mathbf{c}|^{2}\cdot\Bigl(-[\tfrac{1}{2}]_{2}^{k}\partial r_{k-1}
 +\sum_{m=0}^{\infty}[m-\tfrac{1}{2}]_{2}^{k}t_{m}\partial
r_{k+m-2}\\
&+\frac{\hbar}{2}\sum_{m=0}^{k-3}(-1)^{m+1}[-m-\tfrac{3}{2}]_{2}^{k}\partial
r_{m}\partial r_{k-m-3} \Bigr)+\frac{\delta_{k1}}{2\hbar}t_{0}^{2},
\end{align*}
where $\mathbf{c}$ satisfies $c_{1}(\mathcal{X})=\mathbf{c}\omega$.
The degree zero orbifold Gromov-Witten potential is:
\begin{multline*}
\mathcal{D}_{\mathcal{X}}^{0,0}=
\exp\Bigl(\frac{1}{\hbar}\sum_{m=0}^{\infty}r_{m}\langle\langle\tau_{m}|\lambda_{0}\rangle\rangle_{0}
+\frac{1}{\hbar}\sum_{l,m}\frac{1}{2}\langle\langle \tau_{l}\tau_{m}\rangle\rangle_{0}\\
+\frac{1}{\hbar}\sum_{g=0}^{\infty}~
\sum_{k_{1},\cdots,k_{\Sigma_{i}n_{i}}\geq 0}\alpha^{1}_{k_{1}}
\cdots\alpha^{N-1}_{k_{\Sigma_{i}
n_{i}}}\langle\langle\widetilde{\tau}_{k_{1}}\cdots
\widetilde{\tau}_{k_{\Sigma_{i}n_{i}}}|\lambda^{2}_{r_{1}}\rangle\rangle_{g}\Bigr).
\end{multline*}
So applying the operator $\mathcal{L}_{k}$ to
$\mathcal{D}_{\mathcal{X}}^{0,0}$ we obtain the result.
\end{proof}

\begin{rmk}
As explained in \cite{GP}, the formula $w(r,s,t)$ is very
complicated. To get formula of Hurwitz-Hodge integrals, it is not
necessary to write $w(r,s,t)$ down.
\end{rmk}

So the degree zero, genus zero Virasoro constraints for $\mathbb{P}(1,1,N)$ is equivalent to the vanishing of $w(r,s,t)$ in \cite{GP} and $z_{g,k_{1},\cdots,k_{\Sigma_{i}n_{i}}}^{k}(\mathbf{t})$.

Let
\begin{equation}\Label{constant2}
\Gamma_{j,g}=\langle\widetilde{\tau}_{a}|\lambda^{2}_{r_{1}}\rangle_{g}=\int_{\overline{H}_{g}}\psi^{a}\lambda^{2}_{r_{1}},\quad \text{where } a:=\sum_{i=1}^{N-1}n_{i}-\sum_{i=1}^{N-1}\frac{2i}{N}n_{i}-2.
\end{equation}

Again write $\mathbf{\Gamma_{g}}:=(\Gamma_{j,g})_{1\leq j\leq \Sigma_{i}n_{i}}$ as a column vector and let $\mathbf{c_{g}}:=(c_{j,g})_{1\leq j\leq \Sigma_{i}n_{i}}$ be another column vector. Let the index $i$ vary from $1$ to $N-1$ and define a $\Sigma_{i}n_{i}\times\Sigma_{i}n_{i}$ square matrix $A=(a_{st})$ by
\begin{equation}\Label{matrix2}
a_{st}:=\begin{cases}\frac{2i}{N}+a&\text{if
}~ M_{i-1}<s=t\leq M_{i}\,;\\
\frac{2i}{N}&\text{if}~M_{i-1}<t\leq M_{i}~\text{and}~s\neq t\,.\end{cases}
\end{equation}
The matrix can be written as follows:
$$
A=\left[
\begin{array}{cccccccccc}
  \tfrac{2}{N}+a&\cdots&\tfrac{2}{N}&\tfrac{4}{N}&\cdots&\tfrac{4}{N}&\cdots&\tfrac{2(N-1)}{N}&\cdots&\tfrac{2(N-1)}{N} \\
  \tfrac{2}{N}&\cdots&\tfrac{2}{N}&\tfrac{4}{N}&\cdots&\tfrac{4}{N}&\cdots&\tfrac{2(N-1)}{N}&\cdots&\tfrac{2(N-1)}{N} \\
  \cdots&\cdots&\cdots&\cdots&\cdots&\cdots&\cdots&\cdots&\cdots&\cdots\\
  \tfrac{2}{N}&\cdots&\tfrac{2}{N}+a&\tfrac{4}{N}&\cdots&\tfrac{4}{N}&\cdots&\tfrac{2(N-1)}{N}&\cdots&\tfrac{2(N-1)}{N} \\
  \tfrac{2}{N}&\cdots&\tfrac{2}{N}&\tfrac{4}{N}+a&\cdots&\tfrac{4}{N}&\cdots&\tfrac{2(N-1)}{N}&\cdots&\tfrac{2(N-1)}{N} \\
  \cdots&\cdots&\cdots&\cdots&\cdots&\cdots&\cdots&\cdots&\cdots&\cdots \\
  \tfrac{2}{N}&\cdots&\tfrac{2}{N}&\tfrac{4}{N}&\cdots&\tfrac{4}{N}+a&\cdots&\tfrac{2(N-1)}{N}&\cdots&\tfrac{2(N-1)}{N} \\
 \cdots&\cdots&\cdots&\cdots&\cdots&\cdots&\cdots&\cdots&\cdots&\cdots \\
 \tfrac{2}{N}&\cdots&\tfrac{2}{N}&\tfrac{4}{N}&\cdots&\tfrac{4}{N}&\cdots&\tfrac{2(N-1)}{N}+a&\cdots&\tfrac{2(N-1)}{N}\\
  \cdots&\cdots&\cdots&\cdots&\cdots&\cdots&\cdots&\cdots&\cdots&\cdots \\
 \tfrac{2}{N}&\cdots&\tfrac{2}{N}&\tfrac{4}{N}&\cdots&\tfrac{4}{N}&\cdots&\tfrac{2(N-1)}{N}&\cdots&\tfrac{2(N-1)}{N}+a
\end{array}
\right].
$$
It is easy to check that $A$ is nonsingular for $a\neq 0$. Let $\mathbf{A}$ be the matrix obtained from $A$ as follows: for an integer $j$ with $M_{i-1}+1\leq j\leq M_i$ for some $1\leq i\leq N-1$, the $j$-th row of $\mathbf{A}$ is obtained by multiplying the $j$-th row of $A$ by
$$\frac{\frac{1}{2}((\sum_{i=1}^{N-1}n_i)-3)!(\frac{2i}{N}-\frac{1}{2})}{(a+\frac{2i}{N}-\frac{1}{2})!\prod_{i=1}^{N-1}(\frac{2i}{N}-\frac{1}{2})^{n_i}}.$$


The linear system
\begin{equation}\label{linearsystem2}
\mathbf{A}\cdot \mathbf{c_{g}}=\mathbf{\Gamma_{g}}
\end{equation}
has a unique solution which represents $c_{j,g}$ as a linear combination of $\Gamma_{j,g}$'s for $1\leq j\leq \Sigma_{i}n_{i}$. 

For integers $1\leq s\leq N-1$ and $1\leq j\leq\Sigma_{i}n_{i}$, let 
$$(k_{j}+\frac{2s}{N})!=\frac{2s}{N}\cdot(1+\frac{2s}{N})\cdots(k_{j}+\frac{2s}{N}).$$

From the vanishing of $z_{g,k_{1},\cdots,k_{\Sigma_{i}n_{i}}}^{k}(\mathbf{t})$ for $k\geq 1$ and $k_{1},\cdots,k_{\Sigma_{i}n_{i}}\geq 0$, we obtain the following theorem. 

\begin{thm}\Label{surfacevir}
We have
$$
\langle\widetilde{\tau}_{k_{1}}\cdots\widetilde{\tau}_{k_{\Sigma_{i}n_{i}}}\tau_{l_{1}}
\cdots\tau_{l_{n}}|\lambda^{2}_{r_{1}}\rangle_{g}=
\sum_{s=1}^{N-1}\sum_{j=M_{s-1}+1}^{M_{s}}\tfrac{\frac{1}{2}(n+\Sigma_{i}
n_{i}-3)!(k_{j}-\frac{1}{2}+\frac{2s}{N})}
{\prod_{j}(l_{j}-\frac{1}{2})!\prod_{b=1}^{N-1}\prod_{j=M_{b-1}+1}^{M_{b}}(k_{j}-\frac{1}{2}+\frac{2b}{N})!}c_{j,g}.
$$
\end{thm}

\begin{proof}
Again we consider the following recursion given by  $\tfrac{1}{[\frac{1}{2}]_{0}^{k}}\partial t_{l_{1}}\cdots \partial t_{l_{n}}z_{g;k_{1},\cdots,k_{\Sigma_{i}n_{i}}}^{k}(0)=0$:

\begin{multline}\Label{recursion2}
0=-\langle\tau_{k+1}\widetilde{\tau}_{k_{1}}\cdots\widetilde{\tau}_{k_{\Sigma_{i}n_{i}}}\tau_{l_{1}}
\cdots\tau_{l_{n}}|\lambda^{2}_{r_{1}}\rangle_{g}
+\sum_{i=1}^{n}\tfrac{[l_{i}-\frac{1}{2}]^{k}_{0}}{[\frac{1}{2}]_{0}^{k}}\langle\widetilde{\tau}_{k_{1}}\cdots\widetilde{\tau}_{k_{\Sigma_{i}n_{i}}}\tau_{l_{1}}
\cdots\tau_{l_{i}+k}\cdots\tau_{l_{n}}|\lambda^{2}_{r_{1}}\rangle_{g}\\
+\sum_{i=1}^{N-1}\sum_{j=M_{i-1}+1}^{M_{i}}\tfrac{[k_{j}+\frac{2i}{N}-\frac{1}{2}]_{0}^{k}}{[\frac{1}{2}]_{0}^{k}}
\langle\widetilde{\tau}_{k_{1}}\cdots\widetilde{\tau}_{k_{n_{i-1}}}\widetilde{\tau}_{k_{n_{i-1}+1}}\cdots\widetilde{\tau}_{k_{j}+k}\cdots
\widetilde{\tau}_{k_{n_{i}}}\widetilde{\tau}_{k_{n_{i}+1}}\cdots\widetilde{\tau}_{k_{\Sigma_{i}n_{i}}}\tau_{l_{1}}
\cdots\tau_{l_{n}}|\lambda^{2}_{r_{1}}\rangle_{g}.
\end{multline}
Dimension constraints for orbifold Gromov-Witten invariants of $\mathbb{P}(1,1,N)$ gives 
\begin{equation}\label{virdim2}
\frac{1}{2}\Bigl(n+\sum_{i=1}^{N-1}
n_{i}-2\Bigr)
=\sum_{i=1}^{n}\Bigl(l_{i}-\tfrac{1}{2}\Bigr)+\sum_{i=1}^{N-1}\sum_{j=M_{i-1}+1}^{M_{i}}\Bigl(k_{j}-\tfrac{1}{2}+\tfrac{2i}{N}\Bigr)+k.
\end{equation}
So from (\ref{linearsystem2}), (\ref{recursion2}) and
(\ref{virdim2}), using the same method as in the proof of Theorem
\ref{curvevir} we finish the proof.
\end{proof}

\begin{rmk}
The initial values in Theorem \ref{surfacevir} are the following Hurwitz-Hodge integrals
$$\int_{\overline{H}_{g}}\lambda_{r_{1}}^2\psi^a,$$
as in (\ref{constant2}), where $\overline{H}_{g}$ is the Hurwitz scheme parametrizing the
admissible $\mathbb{Z}_{N}$-covers over $\mathbb{P}^{1}$. These integrals can be interpreted as certain orbifold Gromov-Witten invariants of $[\com^2/\mathbb{Z}_N]$. They have been computed in \cite{ccit2}.
\end{rmk}

\section{The degree zero Virasoro conjecture for threefold.}\label{threefolds}

The Virasoro constraints in degree zero for threefolds do not give anything new:
all descendent invariants are reduced to primary ones by string and dilaton equations.
Indeed, the formula for virtual dimension implies that, if there is an insertion of the form $a\psi^k$ with $k\geq 2$,
then there must be another insertion $1$. Thus the exponent in the descendent insertions are reduced by string equation.
If $k=1$, then the class $a$ is either $1$ or not. In the latter case there must be another
insertion $1$ and the string equation is applied again. If $a=1$, then such insertion is removed by the dilaton equation.

In case of $\mathbb{P}(1,1,1,3)$, the relevant primary invariants are related to
some Hurwitz-Hodge integrals arising in orbifold Gromov-Witten theory of $[\com^3/\zz_3]$.
For example, the first one is
\begin{equation*}
\int_{\overline{\mathcal{M}}_{n_{2}+g+2}(\mathcal{B}\mathbb{Z}_{3})}\lambda^{3}_{r_{1}},
\end{equation*}
where $n_{1}+n_{2}=g+2$ and $n_{1},n_{2}$ represent $n_{1}$ stacky points with type $\omega$
and $n_{2}$ stacky points with type $\overline{\omega}$.

These integrals have been predicted in physics \cite{ABK}. For example, let $n_{1}=3$, $n_2=0$, then 
$$\int_{\overline{\mathcal{M}}_{3}(\mathcal{B}\mathbb{Z}_{3})}\lambda^{3}_{r_{1}}=\int_{\overline{H}_{g}(\omega,\omega,\omega)}\lambda^{3}_{r_{1}}=\frac{1}{3}.$$
The generating function of these integrals is computed in \cite{ccit}, \cite{ccit2}.



\appendix

\section{On Conjecture \ref{orbLW}}\label{conj_orbLW}
Define the {\em double inertia orbifold} $II\X$ of $\X$ to be the inertia orbifold of the orbifold $I\X$. A point of $II\X$ is a triple $(x,g,h)$ where $x\in \X, g,h \in Stab_\X(x)$ with $gh=hg$. There is a natural projection $\pi_2:II\X\to \X$ which forgets $g,h$. The genus one invariants participating Conjecture \ref{orbLW} are conjecturally evaluated as follows:

\begin{conjecture}\label{genus_one_inv}
\hfill
\begin{enumerate}
\item
We have $$\<\psi\>_{1,1,0}^\X=\frac{1}{24}\int_{II\X}c_{top}(T_{II\X})=\frac{1}{24}\chi_{top}(IX);$$
\item
For $D\in H^2(\X)$, we have $$\<D\>_{1,1,0}^\X=\frac{1}{24}\int_{II\X}\pi_2^*(D)c_{top-1}(T_{II\X}).$$
\end{enumerate}
\end{conjecture}
This conjecture is formulated based on an analysis of the moduli stack $\Mbar_{1,1}(\X,0)$ of genus one, degree zero orbifold stable maps with one non-stacky marked point. Conjecture (\ref{genus_one_inv}) holds for the case $\X=\mathcal{B}G$, by the work \cite{jk}.

Conjecture (\ref{genus_one_inv}) will be addressed elsewhere.

We now return to Conjecture \ref{orbLW}. Note that $$\frac{1}{4}str\left(\frac{1}{4}-\mu^2\right)=\frac{1}{16}\chi_{top}(IX)-\frac{1}{4}str(\mu^2)$$ because the term $str(1)=\chi_{top}(IX)$, the topological Euler characteristic of the {\em coarse space} of the inertia orbifold $I\X$. Now by Lefschetz trace formula, we have $$\chi_{top}(IX)=\int_{II\X} c_{top}(T_{II\X}).$$ Combining with (\ref{genus_one_inv}), we may rewrite (\ref{orbLW}) as \\\\
{\bf Conjecture \ref{orbLW}':}
$$str(\mu^2)=\frac{1}{12}\int_{II\X}\text{dim}\X c_{top}(T_{II\X})+2c_1(T_\X)c_{top-1}(T_{II\X}).$$

It is clear that $$str(\mu^2)=\sum_{i\in \sI}\sum_{p\geq 0}(-1)^p\left(p+age(\X_i)-\frac{\text{dim}\X}{2}\right)^2\chi(\X_i,\Omega_{\X_i}^p).$$

In case when the orbifold is actually a compact K\"ahler {\em manifold} $X$, (\ref{orbLW}') is reduced to
\begin{equation}\label{LW}
\sum_{p\geq 0}\left(p-\frac{\text{dim}X}{2}\right)^2\chi(X,\Omega_X^p)=\frac{1}{12}\int_X \text{dim}X c_{top}(T_X)+2c_1(T_X)c_{top-1}(T_X),
\end{equation}
 which is a theorem of Libgober and Wood \cite{LW} (see also \cite{G}). (\ref{LW}) is proven by the use of Hirzebruch-Riemann-Roch formula. Presumably (\ref{orbLW}') can be proven by Hirzebruch-Riemann-Roch formula for orbifolds. We plan to address this in the future.

An evidence of Conjecture \ref{orbLW}' is provided by a result of V. Batyrev \cite{ba}, which we now explain.

Let $\X$ be a Gorenstein proper Deligne-Mumford stack with projective coarse moduli space $X$ such that the natural map $\X\to X$ is proper and birational, and the pullback of $K_X$ coincides with $K_\X$. Then a result of T. Yasuda \cite{ya} asserts that Batyrev's {\em stringy Hodge numbers} $h_{st}^{p,q}(X)$ coincide with {\em orbifold Hodge numbers} $h_{orb}^{p,q}(\X):=\text{dim}H_{orb}^{p,q}(\X, \com)$. In terms of generating functions, we have $$E_{st}(X, s,t)=E_{orb}(\X,s,t),$$ where $E_{st}(X,s,t):=\sum_{p, q\geq 0}(-1)^{p+q}h_{st}^{p,q}(X)s^pt^q$ is the stringy E-polynomial and $E_{orb}(\X,s,t):=\sum_{p, q\geq 0}(-1)^{p+q}h_{orb}^{p,q}(\X)s^pt^q$ is the orbifold E-polynomial. Combining this with Corollary 3.10 of \cite{ba} we find that
\begin{equation}\label{evidence}
\sum_{i\in \sI}\sum_{p\geq 0}(-1)^p\left(p+age(\X_i)-\frac{\text{dim}\X}{2}\right)^2\chi(\X_i,\Omega_{\X_i}^p)=\frac{\text{dim}X}{12}e_{st}(X)+\frac{1}{6}c_{st}^{1,n-1}(X),
\end{equation}
where $c_{st}^{1,n-1}(X)$ is defined in \cite{ba}, Definition 3.1.

Since $e_{st}(X)=\chi_{orb}(\X)=\chi_{top}(IX)$, the first terms of the right-hand sides of (\ref{orbLW}') and (\ref{evidence}) coincide. The second terms on the right-hand sides of (\ref{orbLW}') and (\ref{evidence}) take very similar forms--in particular, the number $c_{st}^{1,n-1}(X)$ can be interpreted as a stringy version of $c_1(X)c_{n-1}(X)$. Note that this proves (\ref{orbLW}') in the Calabi-Yau case, because both $c_{st}^{1,n-1}(X)$ and $\int_{II\X}c_1(T_\X)c_{top-1}(T_{II\X})$ vanish when $c_1(T_\X)=0$. 

The following conjecture is of interests in its own right.
\begin{conjecture}
$$c_{st}^{1,n-1}(X)=\int_{II\X}c_1(T_\X)c_{top-1}(T_{II\X}).$$
\end{conjecture}

\subsection{An Example: Weighted Projective Line}
We verify (\ref{orbLW}') for $\X=\mathbb{P}(a,b)$, the weighted projective line with $a,b$ coprime.

First, we have $I\mathbb{P}(a,b)=\mathbb{P}(a,b)\cup\bigcup_{i=1}^{a-1}(B\mathbb{Z}_a)_i\cup\bigcup_{j=1}^{b-1}(B\mathbb{Z}_b)_j$. The component $\mathbb{P}(a,b)$ has age $age(\mathbb{P}(a,b))=0$. Each component $(B\mathbb{Z}_a)_i\simeq B\mathbb{Z}_a$ has age $age((B\mathbb{Z}_a)_i)=i/a$. Each component $(B\mathbb{Z}_b)_j\simeq B\mathbb{Z}_b$ has age $age((B\mathbb{Z}_b)_j)=j/b$. From this we find that the left-hand side of (\ref{orbLW}') is $$(-1/2)^2+(1-1/2)^2+\sum_{i=1}^{a-1}(i/a-1/2)^2+\sum_{j=1}^{b-1}(j/b-1/2)^2$$
$$=\frac{1}{4}+\frac{1}{4}+\frac{(a-1)(2a-1)}{6a}-\frac{a-1}{2}+\frac{a-1}{4}+\frac{(b-1)(2b-1)}{6b}-\frac{b-1}{2}+\frac{b-1}{4}$$
$$=\frac{a+b}{12}+\frac{1}{6}(\frac{1}{a}+\frac{1}{b}).$$

Now we consider the right-hand side of (\ref{orbLW}'). The first term in the right-hand side of (\ref{orbLW}') is $$\frac{1}{12}\int_{II\mathbb{P}(a,b)}c_{top}(T_{II\mathbb{P}(a,b)})=\frac{1}{12}\chi_{top}(\overline{I\mathbb{P}(a,b)})=\frac{a+b}{12}$$ by Gauss-Bonnet.
Now, note that there is a unique component of positive dimension in $II\mathbb{P}(a,b)$, and this component is isomorphic to $\mathbb{P}(a,b)$. Only this component contributes to the second term of the right-hand side of (\ref{orbLW}'). This gives $$\frac{1}{6}\int_{\mathbb{P}(a,b)}c_1(T_{\mathbb{P}(a,b)})=\frac{1}{6}(\frac{1}{a}+\frac{1}{b}).$$
Now it is evident that both sides agree.

\end{document}